\documentclass{amsart}
\usepackage{amssymb}
\usepackage{mathrsfs}
\usepackage{stmaryrd}
\usepackage{bbm}
\usepackage{oldgerm}
\usepackage[english]{babel}
\usepackage[T1]{fontenc}
\usepackage[latin1]{inputenc}
\usepackage[all]{xy}
\usepackage{hyperref}
\usepackage{upref}
\usepackage{xcolor}
\usepackage{tikz-cd}

\hypersetup{
 colorlinks,
 linkcolor={red!50!black},
 citecolor={blue!50!black},
 urlcolor={blue!80!black}
}

\newtheorem{teo}[subsection]{Theorem}
\newtheorem{prop}[subsection]{Proposition}
\newtheorem{cor}[subsection]{Corollary}
\newtheorem{lem}[subsection]{Lemma}

\theoremstyle{definition}

\newtheorem{defi}[subsection]{Definition}
\newtheorem{rema}[subsection]{Remark}

\newtheorem{exemple}[subsection]{Example}

\numberwithin{equation}{subsection}

\newcommand{\Art}{\mathrm{Art}}
\newcommand{\bk}{\overline{k}}
\newcommand{\bta}{\overline{\eta}}
\newcommand{\bs}{\overline{s}}

\newcommand{\cyc}{\mathrm{cyc}}

\newcommand{\et}{\acute{\mathrm{e}}\mathrm{t}}
\newcommand{\Et}{\acute{\mathrm{E}}\mathrm{t}}

\newcommand{\F}{\mathcal{F}}

\newcommand{\G}{\mathcal{G}}

\newcommand{\Gal}{\mathrm{Gal}}
\newcommand{\Hom}{\mathrm{Hom}}
\newcommand{\Ind}{\mathrm{Ind}}
\newcommand{\id}{\mathrm{id}}

\newcommand{\Lc}{\mathcal{L}}

\newcommand{\Mod}{\mathrm{Mod}}

\newcommand{\m}{\mathfrak{m}}

\newcommand{\Ow}{\mathcal{O}}

\newcommand{\Pic}{\mathrm{Pic}}

\newcommand{\Res}{\mathrm{Res}}
\newcommand{\rk}{\mathrm{rk}}
\newcommand{\Sh}{\mathrm{Sh}}

\newcommand{\Ver}{\mathrm{Ver}}
\newcommand{\ver}{\mathrm{ver}}

\title{Geometric local $\varepsilon$-factors in higher dimensions}
\author{Quentin Guignard}
\address{Institut des
Hautes \'Etudes Scientifiques, 35 route de Chartres, 91440 Bures-sur-Yvette, France}
\address{
\'Ecole Normale Sup\'erieure, 45 rue d'Ulm, 75005 Paris, France}
\email{quentin.guignard@ens.fr}

\begin{document}

\begin{abstract} We prove a product formula for the determinant of the cohomology of an $\ell$-adic sheaf over an arbitrary proper scheme over a perfect field of positive characteristic $p$ distinct from $\ell$. The local contributions are constructed by iterating vanishing cycle functors as well as certain ``refined Artin conductors'', the latter being exact additive functors which can be considered as linearized versions of Artin conductors and local $\varepsilon$-factors. We provide several applications of our higher dimensional product formula, such as twist formulas for global $\varepsilon$-factors.

\end{abstract}

\maketitle
\tableofcontents


\section{Introduction \label{intro}}

\subsection{\label{0.0.1}} Let $k$ be a perfect field of positive characteristic $p$ and let $\bk$ be an algebraic closure of $k$. Let $\Lambda$ be either $\overline{\mathbb{F}}_{\ell}$ or $\overline{\mathbb{Q}}_{\ell}$ for a fixed prime number $\ell$ distinct from $p$. The \textit{$\varepsilon$-factor} associated to a pair $(X,\F)$, where $X$ is a proper $k$-scheme and $\F$ is a bounded constructible complex of \'etale sheaves of $\Lambda$-modules on $X$, is the Galois line given by
$$
\varepsilon_{\bk}(X,\F) = \det(R\Gamma(X_{\bk} ,\F))^{-1}.
$$
In this paper, we simply consider this object as a homomorphism from $G_k = \Gal(\bk/k)$ to $\Lambda^{\times}$. We will also consider the Euler characteristic
$$
\chi(X,\F) = \rk(R\Gamma(X_{\bk} ,\F)).
$$
In the case where $X$ is a smooth curve over a finite field, it was conjectured by Langlands that the global $\varepsilon$-factor $\varepsilon_{\bk}(X,\F) $ should split as a product of local contributions. This conjecture was motivated by the corresponding theory of automorphic local $\varepsilon$-factors. Partial results were obtained by Dwork \cite{Dw} and Deligne \cite{De73}, and the product formula for $\ell$-adic sheaves on a smooth curve over a finite field was proved in its full generality by Laumon (\cite{La87}, Th. 3.2.1.1). The case of a smooth curve over an arbitrary perfect field of positive characteristic was handled in \cite{G19}. An alternative approach by ``spreading out'' in order to reduce to the finite field case was also given by Yasuda in \cite{Ya1}, \cite{Ya2}, \cite{Ya3}.

In the geometric setting of \cite{G19}, the local $\varepsilon$-factors are realized as determinants of certain Galois modules arising in the cohomology of Gabber-Katz extensions, cf. (\cite{G19}, 9.2). These Galois modules are thus \textit{linearizations} of local $\varepsilon$-factors: these are given by additive functors, whose determinants yield the local factors. In this paper we show how such linearizations, henceforth labelled \emph{refined Artin conductors}, allow to split higher dimensional global $\varepsilon$-factors into a finite product of local contributions. The resulting higher local $\varepsilon$-factors arise as iteration of vanishing cycle functors and of these refined Artin conductors. Our main result is the following:

\begin{teo}\label{teo0}Let $k$ be a perfect field of positive characteristic, and let $X$ be a proper $k$-scheme. Let $\Lambda$ be either $\overline{\mathbb{F}}_{\ell}$ or $\overline{\mathbb{Q}}_{\ell}$, where $\ell$ is a prime number invertible in $k$. Then there exists a collection
$$
(E(X_{(x)}, - ) )_{x \in |X|},
$$ 
of triangulated functors $E(X_{(x)}, - ) : D_c^b(X_{(x)}, \Lambda) \rightarrow D_c^b(x,\Lambda)$, indexed by the set of closed points of $X$, where we denoted by $X_{(x)}$ the henselization of $X$ at a point $x$, such that:
\begin{enumerate}
\item for any object $\F$ of $D_c^b(X,\Lambda)$, we have $E(X_{(x)}, \F_{|X_{(x)}} ) \simeq 0$ for all but finitely many closed points $x$ of $X$;
\item for any closed point $x$ of $X$, any object $\F$ of $D_c^b(X_{(x)},\Lambda)$, and any object $\G$ of $D_c^b(x,\Lambda)$, we have an isomorphism
$$
E(X_{(x)}, \F \otimes \mathrm{sp}^{-1} \G) \simeq E(X_{(x)}, \F) \otimes \G,
$$
where $\mathrm{sp} : X_{(x)} \rightarrow x$ is the canonical specialization morphism, and this isomorphism is functorial in $\F$ and $\G$;
\item for any object $\F$ of $D_c^b(X,\Lambda)$, we have
\begin{align*}
\varepsilon_{\bk}(X,\F) &= \det \left( \bigoplus_{x \in |X|} \Ind_{G_x}^{G_k} E(X_{(x)}, \F_{|X_{(x)}} ) \right), \\
-\chi(X,\F) &= \rk \left( \bigoplus_{x \in |X|} \Ind_{G_x}^{G_k} E(X_{(x)}, \F_{|X_{(x)}} ) \right),
\end{align*}
where $G_k = \Gal(\bk/k)$ is the Galois group of an algebraic closure of $k$ and $G_x$ is the Galois group of $\bk / k(x)$.
\end{enumerate}
\end{teo}

This theorem will be proved in Section \ref{prf}. It will be clear from the proof that the resulting functors $E(X_{(x)}, - )$ are not uniquely determined by the conditions in the conclusion of Theorem \ref{teo0}. The factorisation of $\varepsilon_{\bk}(X,\F)$ in Theorem \ref{teo0} can be written as
$$
\varepsilon_{\bk}(X,\F) = \prod_{x \in |X|} \delta_{x/k}^{\rk  (E(X_{(x)}, \F_{|X_{(x)}} ))}  \det( E(X_{(x)}, \F_{|X_{(x)}} ) \circ \mathrm{ver}_{x/k},
$$
where $\delta_{x/k} : G_k \rightarrow \{ \pm1 \}$ is the signature homomorphism for the left action of $G_k$ on $G_k/ G_{k(x)}$, and $\mathrm{ver}_{x/k} : G_x^{\mathrm{ab}} \rightarrow G_k^{\mathrm{ab}}$ is the transfer homomorphism associated to the finite index subgroup $G_x$ of $G_k$. The factors 
$$
 \det \left( E(X_{(x)}, \F_{|X_{(x)}} ) \right),
$$
thus play the role of local $\varepsilon$-factors in higher dimension. 

\subsection{\label{0.0.2}} Theorem \ref{teo0} admits the following immediate consequence:


\begin{teo} \label{teo0.1} Let $X$ be a proper $k$-scheme. Let $\F_1$ and $\F_2$ be objects of $D_c^b(X,\Lambda)$ such that $\F_{1 | X_{(x)}}$ is isomorphic to $\F_{2 | X_{(x)}}$ for any closed point $x$ of $X$. Then we have
\begin{align*}
\chi(X,\F_1) &= \chi(X,\F_2), \\
\varepsilon_{\bk}(X,\F_1) &= \varepsilon_{\bk}(X,\F_2).
\end{align*}
\end{teo}

 The conclusion of \ref{teo0.1} regarding Euler characteristics is originally due to Deligne, whose proof was written by Illusie in \cite{Ill81}. The part concerning $\varepsilon$-factors appears to be new, though the case of a finite base field alternatively follows from the Grothendieck-Lefschetz trace formula, as was pointed out to the author by Takeshi Saito.

\begin{teo} \label{teo0.2} Let $S$ be a henselian trait of equicharacteristic $p$, with closed point $s$ such that $k(s)$ is perfect. Let $f : X \rightarrow S$ be a proper morphism and let $\F_1$ and $\F_2$ be objects of $D_c^b(X,\Lambda)$ such that $\F_{1 | X_{(x)}}$ is isomorphic to $\F_{2 | X_{(x)}}$ for any closed point $x$ of the special fiber $X_s$. Then we have
\begin{align*}
a(S,Rf_* \F_1) &= a(S,Rf_*\F_2), \\
\varepsilon_{\overline{s}}(S,Rf_* \F_1, \omega) &= \varepsilon_{\overline{s}}(S,Rf_* \F_2, \omega).
\end{align*}
for any non zero meromorphic $1$-form $\omega$ on $S$. Here we denoted by $a(S,-)$ the Artin conductor, cf. $($\cite{G19}, $7.2)$, and by $\varepsilon_{\overline{s}}(S,-,\omega)$ the geometric local $\varepsilon$-factor, cf. $($\cite{G19}, $9.2)$.
\end{teo}

This a consequence of Theorem \ref{teo0.2prime}, whose proof relies on Theorem \ref{teo0}.


\begin{teo}[Twist Formula] \label{teo0.3} Let $X$ be a proper $k$-scheme and let $(E(X_{(x)},-)_{x \in |X|}$ be as in Theorem \ref{teo0}. Let $\F$ be an object of $D_c^b(X,\Lambda)$, and let $\G$ be a $\Lambda$-local system of constant rank $r$ on $X$. We then have
\begin{align*}
\chi(X,\F \otimes \G) &= r \chi(X,\F) \\
\varepsilon_{\bk}(X,\F \otimes \G) &= \varepsilon_{\bk}(X,\F)^r \prod_{x \in |X|} \left( \det(\G_{\overline{x}}) \circ \ver_{x/k} \right)^{\rk(E(X_{(x)},\F))}
\end{align*}
\end{teo}

This is an immediate consequence of Theorem \ref{teo0}. A similar twist formula was obtained in \cite{UYZ} in the case of a projective smooth scheme over a finite field. Let us set
$$
cc_X(\F) = - \sum_{x \in |X|} \rk(E(X_{(x)},\F)) [x],
$$
as a $0$-cycle on $X$. Theorem \ref{teo0} yields $\deg(cc_X(\F)) = \chi(X,\F)$, and the conclusion of Theorem \ref{teo0.3} can be written as 
\begin{align*}
\varepsilon_{\bk}(X,\F \otimes \G) = \varepsilon_{\bk}(X,\F)^r \langle \det(\G) , -cc_X(\F) \rangle,
\end{align*}
where the pairing $\langle -, - \rangle $ is defined by 
\begin{align}\label{pairing}
\langle \Lc , c \rangle = \prod_i (\Lc_{\overline{x_i}} \circ \ver_{x_i/k})^{n_i},
\end{align}
for any $\Lambda$-local system $\Lc$ of rank $1$ on $X$ and any $0$-cycle $c=\sum_i n_i [x_i]$. When $k$ is finite and when $X$ is projective smooth, then the main result of \cite{UYZ} and the class field theory of Kato-Saito \cite{KS} yield that $cc_X(\F)$ coincides in the Chow group $CH_0(X)$ with the characteristic class defined in (\cite{Sai}, Def. 5.7). When $k$ is arbitrary and $X$ is projective smooth, one can simply assert that the difference between $cc_X(\F)$ and Saito's characteristic class belong to the kernel of the pairing \ref{pairing}. As in (\cite{UYZ}, 5.26), this implies that the formation of Saito's characteristic class on projective smooth $k$-schemes commutes with proper pushforward, up to a $0$-cycle in the kernel of the pairing \ref{pairing}. Finally, we notice that Theorem \ref{teo0} actually yields a slightly stronger twist formula:

\begin{teo} \label{teo0.4}Let $X$ be a proper $k$-scheme and let $(E(X_{(x)},-)_{x \in |X|}$ be as in Theorem \ref{teo0}. Let $r$ be an integer. Let $\F_1$ and $\F_2$ be objects of $D_c^b(X,\Lambda)$ such that for any closed point $x$ of $X$, we have $\F_{2|X_{(x)}} \simeq \F_{1|X_{(x)}} \otimes \G_x$ for some free $\Lambda$-module $\G_x$ of rank $r$ endowed with an admissible action of $\Gal(\bk/k(x))$. We then have
\begin{align*}
\chi(X,\F_2) &= r \chi(X, \F_1) \\
\varepsilon_{\bk}(X,\F_2) &= \varepsilon_{\bk}(X,\F_1)^r \prod_{x \in |X|} \left( \det(\G_{x}) \circ \ver_{x/k} \right)^{\rk(E(X_{(x)},\F))}.
\end{align*}
\end{teo}

\subsection{\label{0.0.5}} Theorem \ref{teo0}, as well as its consequences \ref{teo0.1}, \ref{teo0.2}, \ref{teo0.3}, \ref{teo0.4}, can be generalized to (complexes of) \'etale $\ell$-adic sheaves twisted by a $2$-cocycle on $G_k$, cf. (\cite{G19}, Sect. 3). We choose to restrict the exposition to non twisted sheaves for the sake of simplicity and concision. Let us however notice that in the twist formula \ref{teo0.3}, we can allow $\G$ to be a twisted local system, and this ultimately yields a more precise conclusion in the discussion following Theorem \ref{teo0.3}, regarding the commutation of the formation of Saito's characteristic class with proper pushforward.

\subsection{\label{0.0.3}} Let us now describe the organization of this paper. In Section \ref{5}, we give a few sorites on Galois equivariant \'etale sheaves on schemes. We also study transformations of such sheaves afforded by \textit{$\Lambda$-linear contractions}, cf. \ref{contract}.

In Section \ref{7}, we give a brief review of nearby and vanishing cycles functors, as presented in (\cite{SGA7}, XIII), and we state a few results regarding the composition of these functors with a $\Lambda$-linear contraction.

In Section \ref{0}, we review the theory of Gabber-Katz extensions, and we use it to provide three families of $\Lambda$-linear contractions: Katz's cohomological construction of the Swan module from (\cite{katz}, 1.6), cf. \ref{3.2}, a functor refining geometric class field theory, cf. \ref{3.3}, and a refined Artin conductor, cf. \ref{0.1}, which linearizes Artin conductors and geometric local $\varepsilon$-factors. Only the second and the third of these constructions will be used in the proof of our main results, though Katz's construction was inspirational to us.

Finally, we prove the main Theorem \ref{teo0} in Section \ref{prf}. We first provide a stronger result in the case of a projective line in \ref{4.0}, by using the product formula from \cite{G19}. We then prove in \ref{4.1} and \ref{4.2} using Chow's lemma that it is sufficient to prove Theorem \ref{teo0} in the case of a projective space. The latter case is handled in \ref{4.3} by choosing an arbitrary pencil and by applying to the latter the product formula from \ref{4.0}, as well as the functoriality properties from Section \ref{7}. 

\subsection{\label{0.0.4} Acknowledgements} This work was prepared at the Institut des Hautes \'Etudes Scientifiques and the \'Ecole Normale Sup\'erieure while the author benefited from their hospitality and support. This text is a continuation of the part \cite{G19} of the author's doctoral thesis, realized under the supervision of Ahmed Abbes. The author would like to thank Ahmed Abbes, Lei Fu, Ofer Gabber, Dennis Gaitsgory, Adriano Marmora, Deepam Patel, Takeshi Saito, Daichi Takeuchi, Enlin Yang for useful discussions.

\subsection{\label{conv}Conventions and notation}
%

 Throughout this paper, we fix distinct prime numbers $p,\ell$. We fix a non trivial homomorphism $\psi : \mathbb{F}_p \rightarrow \overline{\mathbb{Z}_{\ell}}^{\times}$, and every $\ell$-adic coefficient ring $\Lambda$ will be assumed to contain the image of $\psi$. We denote by $\Lc_{\psi}$ the corresponding Artin-Schreier sheaf on $\mathbb{A}^1_{\mathbb{F}_p}$. For any $\mathbb{F}_p$-scheme $X$ and any global function $f$ in $\Gamma(X,\Ow_X)$, we denote by $\Lc_{\psi} \lbrace f \rbrace$ the pullback of $\Lc_{\psi}$ by the morphism $f : X \rightarrow \mathbb{A}^1_{\mathbb{F}_p}$.

We write ``qcqs'' for ``quasi-compact quasi-separated'', and for any scheme $S$ we denote by $\Et_S$ the category of qcqs \'etale $S$-schemes, and by $S_{\et}$ the \'etale topos of $S$, namely the topos of sheaves of sets on the site $\Et_S$. 

\section{Preliminaries \label{5}}

In this section, we record various sorites regarding equivariant sheaves and transformations thereof. Throughout this section, we let $s$ be the spectrum of a field and we let $\bs$ be a separable closure of $s$. We denote by $G_s$ the Galois group of $k(\bs)$ over $k(s)$. In particular, the group $G_s$ acts on the right on $\bs$.

\subsection{\label{5.1}} Let $\theta : Q \rightarrow G_s$ be a continuous homomorphism of profinite groups and let $Y$ be an $s$-scheme.

\begin{defi}[\cite{SGA7} XIII 1.1]\label{equivsheaf} Let $\Lambda$ be a ring. Let $\F$ be an \'etale sheaf of sets (resp. of $\Lambda$-modules) on $Y_{\bs}$. 
\begin{enumerate}

\item A (left) \textit{action of $Q$ compatible with} $\theta$ on $\F$ is a collection $(\sigma(q))_{g \in Q}$ of isomorphisms $\sigma(q) : \F \rightarrow \theta(q)_*\F$, for each $q$ in $Q$, such that 
$$
\theta(q_2)_* \sigma(q_1)  \circ \sigma(q_2) = \sigma(q_1 q_2),
$$
for any $q_1,q_2$ in $Q$. 
\item An action of $Q$ compatible with $\theta$ on $\F$ is \textit{continuous} if for any qcqs \'etale $Y$-scheme $U$, the resulting left action of $Q$ on $\F(U_{\bs})$ is continuous, when the latter is endowed with the discrete topology. 
\item We denote by $\Sh(Y,Q,\theta)$ (resp. $\Sh(Y,Q,\theta, \Lambda)$), or simply by $\Sh(Y,Q)$ (resp. $\Sh(Y,Q, \Lambda)$) if no confusion can arise from this abuse of notation, the category of \'etale sheaves (resp. \'etale sheaves of $\Lambda$-modules) on $Y_{\bs}$ endowed with a continuous action of $Q$ compatible with $\theta$. The morphisms in this category are the morphisms of sheaves on $Y_{\bs}$ which commute with the action of $Q$.
\end{enumerate}
\end{defi}

If one takes $Y$ to be the base scheme $s$, then the category $\Sh(Y,Q,\theta)$ (resp. $\Sh(Y,Q,\theta, \Lambda)$) is simply the category of sets (resp. of $\Lambda$-modules) endowed with a continuous left action of $Q$, which we also denoted by $Q-\mathrm{Sets}$ (resp. $Q-\mathrm{Mod}_{\Lambda}$). 

\begin{rema}\label{rema5.1.0} One can alternatively describe $\Sh(Y,Q,\theta)$ (resp. $\Sh(Y,Q,\theta, \Lambda)$) as the product topos $\mathrm{B}Q \times_{s_{\et}} Y_{\et}$ (resp. as the category of $\Lambda$-modules in the topos $\mathrm{B}Q \times_{s_{\et}} Y_{\et}$), where $\mathrm{B}Q$ is the topos of sets endowed with a continuous left action of $Q$. In particular, the abelian category $\Sh(Y,Q,\theta, \Lambda)$ has enough injectives.
\end{rema}

\begin{prop} The natural functor $\Sh(Y) \rightarrow \Sh(Y,G_s,\mathrm{id})$ is an equivalence of categories.
\end{prop}

This is (\cite{SGA7} XIII 1.1.3(ii)). A quasi-inverse to this functor can be described as follows: to any object $\F$ of $\Sh(Y,G_s,\mathrm{id})$ one associates the \'etale sheaf on $Y$ whose sections on a quasi-compact \'etale $Y$-scheme $U$ are given by the $G_s$-invariants $\F(U_{\bs})^{G_s}$.

\subsection{\label{5.5}}Let $\Lambda$ be a ring. Let $Y$ be an $s$-scheme and let $\bs \rightarrow s' \rightarrow s$ be a factorization of $\bs \rightarrow s$, where $s'$ is a finite extension of $s$.  Let us consider a commutative diagram

\begin{center}
 \begin{tikzpicture}[scale=1]

\node (A) at (0,0) {$Q'$};
\node (B) at (2,0) {$G_{s'}$};
\node (C) at (0,1) {$Q$};
\node (D) at (2,1) {$G_s$};

\path[->,font=\scriptsize]
(A) edge node[above]{$\theta_{|Q'}$} (B)
(A) edge  (C)
(B) edge (D)
(C) edge node[above]{$\theta$} (D);
\end{tikzpicture} 
\end{center}
of continuous homomorphisms of profinite groups, with injective vertical arrows. We then have a restriction functor
$$
\Res_{Q'}^Q : \Sh(Y,Q, \theta, \Lambda) \rightarrow \Sh(Y,Q',\theta_{|Q'}, \Lambda),
$$
obtained by restricting to $Q'$ the action of $Q$ on a given object of $\Sh(Y,Q, \Lambda)$. This is an exact and faithful functor. If $Q'$ is an open subgroup of $Q$, then the function $\Res_{Q'}^Q$ admits a left adjoint, the induction functor
$$
\Ind_{Q'}^Q :  \Sh(Y,Q',\theta_{|Q'}, \Lambda) \rightarrow \Sh(Y,Q,\theta, \Lambda),
$$
which is exact as well. 

\begin{exemple}\label{generator} Let $U$ be an object of $\Et_{Y_{s'}}$. Then 
$$
\Lambda_{Q',U} = \Ind_{Q'}^Q \Lambda [ \Hom_{Y_{\bs}}(-,U \times_{s'} \bs)],
$$
is an object of $\Sh(Y,Q,\Lambda)$. For any object $\F$ of $\Sh(Y,Q,\Lambda)$, the natural homomorphism
$$
\Hom(\Lambda_{U,Q'}, \F) \rightarrow \F(U \times_{s'} \bs)^{Q'},
$$
is an isomorphism.
\end{exemple}

\begin{prop}\label{acyclicity1} Let $\F$ be an injective object in $\Sh(Y,Q, \Lambda)$. Let $\mathcal{U} = (U_i \rightarrow U)_{i \in I}$ be a finite cover in $\Et_{Y_{\bs}}$. Then for any $j > 0$, the \v{C}ech cohomology group $H^j(\mathcal{U},\F)$ vanishes.
\end{prop}

Since $I$ is finite, and since any object of $\Et_{Y_{\bs}}$ is qcqs, there exists a factorization $\bs \rightarrow s' \rightarrow s$ as in \ref{5.5} and a cover $\mathcal{U}' = (U_i' \rightarrow U')_{i \in I}$ in $\Et_{Y_{s'}}$ such that $\mathcal{U}$ is isomorphic to $\mathcal{U}' \times_{s'} \bs = (U_i' \times_{s'} \bs \rightarrow U' \times_{s'} \bs)_{i \in I}$. We henceforth assume that $\mathcal{U}$ is actually equal to $\mathcal{U}' \times_{s'} \bs$.

For any open subgroup $Q'$ of $Q$ whose image by $\theta$ is contained in $G_{s'}$, let us consider the chain complex
$$
\Lambda_{Q',\mathcal{U}'} : \cdots \rightarrow \bigoplus_{i_0,i_1 \in I} \Lambda_{Q',U_{i_0}' \times_{U'} U_{i_1}'}  \rightarrow \bigoplus_{i_0\in I} \Lambda_{Q', U_{i_0}'} \rightarrow 0	\rightarrow \cdots,
$$
where the last nonzero term is placed in degree $0$, cf. Example \ref{generator} for the notation. If $\check{\mathcal{C}}^\bullet (\mathcal{U}, \F)$ is the \v{C}ech complex of $\F$ with respect to the cover $\mathcal{U}$, then the natural homomorphism
$$
\Hom(\Lambda_{Q',\mathcal{U}'},\F) \rightarrow \check{\mathcal{C}}^\bullet (\mathcal{U}, \F)^{Q'},
$$ 
is an isomorphism, cf. Example \ref{generator}. The natural homomorphism $\Lambda_{Q',U'}[0] \rightarrow \Lambda_{Q',\mathcal{U}'}$ is a quasi-isomorphism by (\cite{stacks-project} 03AT) and by exactness of the induction functor $\Ind_{Q'}^Q$. By applying the exact functor $\Hom(-,\F)$, this yields that the natural homomorphism from $\F(U)^{Q'}$ to $\check{\mathcal{C}}^\bullet (\mathcal{U}, \F)^{Q'}$ is a quasi-isomorphism. The conclusion then follows by taking the colimit over all open subgroups $Q'$ of $Q$ with image by $\theta$ contained in $G_{s'}$.

\subsection{\label{5.2}} Let $\theta : Q \rightarrow G_s$ be a continuous homomorphism of profinite groups and let $\Lambda$ be a ring. Let $f : Y' \rightarrow Y$ be a quasi-compact morphism of $s$-schemes. For any object $\F$ of $\Sh(Y,Q,\theta, \Lambda)$, with corresponding action $(\sigma(q))_{q \in Q}$, the data $(f_*\sigma(q))_{q \in Q}$ is a continous action of $Q$ on $f_*\F$  compatible with $\theta$. We thus obtain a functor
$$
f_* : \Sh(Y',Q,\theta, \Lambda) \rightarrow \Sh(Y,Q,\theta, \Lambda).
$$
We similarly have a functor
$$
f^{-1} : \Sh(Y,Q,\theta, \Lambda) \rightarrow \Sh(Y',Q,\theta, \Lambda),
$$
which is a left adjoint to $f_*$.

\begin{rema}\label{rema5.2.0} With the description from Remark \ref{rema5.1.0}, this pair of adjoint functors is induced by the morphism of toposes $\mathrm{id} \times f$ from $\mathrm{B}Q \times_{s_{\et}} Y'_{\et}$ to $\mathrm{B}Q \times_{s_{\et}} Y_{\et} $.
\end{rema}

\subsection{\label{5.3}} Let $\theta : Q \rightarrow G$ be a surjective continuous homomorphism of profinite groups and let $I$ be the its kernel. Let $\Lambda$ be a ring. For any element $q$ of $Q$ and any object $M$ of $I-\Mod_{\Lambda}$, we denote by $q^* M$ the object of $I-\Mod_{\Lambda}$ with underlying $\Lambda$-module $M$, on which an element $t$ of $I$ acts as $m \mapsto qtq^{-1} m$. This yields a $\Lambda$-linear functor
$$
q^* : I-\Mod_{\Lambda} \rightarrow I-\Mod_{\Lambda},
$$ 
and we have $(q_1 q_2)^* = q_2^* q_1^*$. For any element $q$ of $I$, we have a $\Lambda$-linear natural transformation
$$
\lambda_q : \id \rightarrow q^*, 
$$
given on an object $M$ of $I-\Mod_{\Lambda}$ by the homomorphism $m \mapsto qm$.

\begin{defi}\label{contract} A\textit{ $\Lambda$-linear $(Q,I)$-contraction} is a $\Lambda$-linear functor
$$
T : I-\Mod_{\Lambda} \rightarrow \Mod_{\Lambda},
$$ 
endowed with a collection $\tau = (\tau(q))_{q \in Q}$ of $\Lambda$-linear natural transformations
$$
\tau(q) : T q^* \rightarrow T,
$$
satisfying the following conditions:
\begin{enumerate}
\item \label{contract1} the functor $T$ is exact and commutes with small filtered colimits;
\item if an object $M$ of $ I-\Mod_{\Lambda} $ is finitely generated as a $\Lambda$-module, then $T(M)$ is a finitely generated $\Lambda$-module, and there exists an open subgroup $Q'$ of $Q$ such that for any element $q$ of $Q'$, we have $q^* M = M$ and $\tau(q)_M = \id_{T(M)}$;
\item for any elements $q_1,q_2$ of $Q$, then the natural transformation $\tau(q_1) \circ \tau(q_2)q_1^*$ from $T q_2^* q_1^*$ to $T$ coincides with $\tau(q_1 q_2)$;
\item for any element $q$ of $I$, we have $\tau(q) \circ T(\lambda_q) = \id_T$. 
\end{enumerate}
\end{defi}

We will give non trivial examples of $\Lambda$-linear $(Q,I)$-contractions in the next section, cf. \ref{3.2}, \ref{3.3} and \ref{0.1}. 

\subsection{\label{5.6}} Let $\theta : Q \rightarrow G$ be a surjective continuous homomorphism of profinite groups with kernel $I$ and let $\Lambda$ be a ring. Let $(T,\tau)$ be a $\Lambda$-linear $(Q,I)$-contraction, cf. \ref{contract}. Let $G'$ be a closed subgroup of $G$ and let $Q'$ be its inverse image by $\theta$. Then one can construct a functor 
$$
T_{G'} : Q'-\Mod_{\Lambda} \rightarrow G'-\Mod_{\Lambda},
$$
as follows. For any element $q$ of $Q'$,  we have a $\Lambda$-linear natural transformation
$$
\lambda_q : \id \rightarrow q^*, 
$$
of endofunctors of $Q'-\Mod_{\Lambda}$, cf. \ref{5.3}. We then obtain a $\Lambda$-linear natural transformation
$$
\tau(q) \circ T(\lambda_q) : T \phi \rightarrow T \phi.
$$
where $\phi$ is the forgetful functor from $Q'-\Mod_{\Lambda}$ to $I-\Mod_{\Lambda}$. Lemma \ref{sanity} below then ensures that for any object $M$ of $Q'-\Mod_{\Lambda}$, one obtains a continuous left action of $Q'$ on $T(M)$ by letting an element $q$ of $Q'$ act as $(\tau(q) \circ T(\lambda_q))_M$. By the axiom \ref{contract} (4), this action factors through $\theta$ and therefore yields an object of $G'-\Mod_{\Lambda}$.

\begin{lem}\label{sanity} For any elements $q_1,q_2$ of $Q'$, we have
$$
\tau(q_1) \circ T(\lambda_{q_1}) \circ \tau(q_2) \circ T(\lambda_{q_2}) = \tau(q_1q_2) \circ T(\lambda_{q_1q_2}).
$$
Moreover, for any object $M$ of $Q'-\Mod_{\Lambda}$ and element $m$ of $T(M)$, there exists an open subgroup $Q''$ of $Q'$ such that $\tau(q)  T(\lambda_q) m = m$ for any element $q$ of $Q''$.
\end{lem}

The first assertion follows from the identities
\begin{align*}
T(\lambda_{q_1}) \circ \tau(q_2) &= \tau(q_2) q_1^* \circ T(q_2^* \lambda_{q_1}), \\
\tau(q_1) \circ \tau(q_2)q_1^* &= \tau(q_1 q_2), \\
q_2^* \lambda_{q_1} \circ \lambda_{q_2} &= \lambda_{q_1 q_2}.
\end{align*}
The first equality stems from the fact that $\tau(q_2)$ is a natural transformation, the second follows from \ref{contract} (3), and the third identity is tautological. For the last assertion of Lemma \ref{sanity}, we can assume by \ref{contract} (1) that $M$ is finitely generated as a $\Lambda$-module, and the conclusion then follows from \ref{contract} (2).

\subsection{\label{5.8}} Since $T$ and $T_G$ are exact functors, we still denote by $T$ and $T_G$ the functors
\begin{align*}
T &: D(I-\Mod_{\Lambda}) \rightarrow D(\Mod_{\Lambda}), \\
T_G &:  D(Q-\Mod_{\Lambda}) \rightarrow D(G-\Mod_{\Lambda}),
\end{align*}
induced by them on the corresponding derived categories.

\begin{prop}\label{tensor} We have a bifunctorial isomorphism
$$
T(M \stackrel{L}{\otimes}_{\Lambda} N) \cong T(M) \stackrel{L}{\otimes}_{\Lambda} N
$$
with $M$ in $D^-(I-\Mod_{\Lambda})$ and $N$ in $D^-(\Mod_{\Lambda})$. We similarly have a bifunctorial isomorphism
$$
T_G(M \stackrel{L}{\otimes}_{\Lambda} N) \cong T_G(M) \stackrel{L}{\otimes}_{\Lambda} N
$$
with $M$ in $D^-(Q-\Mod_{\Lambda})$ and $N$ in $D^-(G-\Mod_{\Lambda})$.
\end{prop}

We only prove the first assertion, the second one being similar. Let $\mathcal{C}$ be the $\Lambda$-linear abelian category of bounded above chain complexes with terms in the full subcategory of $\Mod_{\Lambda}$ consisting of $\Lambda$-modules of the form $\bigoplus_{i < \kappa} \Lambda$ for some cardinal $\kappa$. Since $T$ is $\Lambda$-linear and commutes with small direct sums by \ref{contract} (1), we have a functorial isomorphism
$$
T(M \otimes_{\Lambda} N) \cong T(M) \otimes_{\Lambda} N
$$
with $M$ in $D^-(I-\Mod_{\Lambda})$ and $N$ in $\mathcal{C}$. This isomorphism factors through the quotient $D^-(\Mod_{\Lambda})$ of $\mathcal{C}$ since for any homomorphism $\lambda : N_1 \rightarrow N_2$ in $\mathcal{C}$, which is homotopic to $0$, the corresponding homomorphisms
\begin{align*}
T(M \otimes_{\Lambda} N_1) &\rightarrow T(M \otimes_{\Lambda} N_2) \\
T(M) \otimes_{\Lambda} N_1 &\rightarrow T(M) \otimes_{\Lambda} N_2,
\end{align*}
vanish in $D^-(\Mod_{\Lambda})$.

\begin{cor}\label{flat} Let $M$ be an object of $I-\Mod_{\Lambda}$ which is flat as a $\Lambda$-module. Then $T(M)$ is a flat $\Lambda$-module.
\end{cor}

Indeed, for any $\Lambda$-module $N$, the complex $M \stackrel{L}{\otimes}_{\Lambda} N$ is quasi-isomorphic to $M \otimes_{\Lambda} N[0]$ by flatness of $M$, hence Proposition \ref{tensor} yields an isomorphism
$$
T(M)  \stackrel{L}{\otimes}_{\Lambda} N \cong T(M \otimes_{\Lambda} N)[0]
$$ 
in $D^-(\Mod_{\Lambda})$. In particular, the complex $T(M)  \stackrel{L}{\otimes}_{\Lambda} N$ is acyclic in positive degrees, hence the flatness of $M$.

\subsection{\label{5.7}} Let $\theta : Q \rightarrow G_s$ be a surjective continuous homomorphism of profinite groups with kernel $I$ and let $\Lambda$ be a ring. Let us consider a $\Lambda$-linear $(Q,I)$-contraction $(T,\tau)$, cf. \ref{contract}. Let $Y$ be an $s$-scheme.

\begin{defi} For any object $\F$ of $\Sh(Y,Q, \Lambda)$, we define $T(\F)$ to be the \'etale sheaf of $\Lambda$-modules on $Y$ whose sections on a qcqs \'etale $Y$-scheme $U$ are given by the $G_s$-invariants $T_{G_s}(\F(U_{\bs}))^{G_s}$ (cf. \ref{5.6}).
\end{defi}

The mentioned presheaf is indeed a sheaf by left exactness of the functor $M \mapsto T_{G_s}(M)^{G_s}$ from $Q-\mathrm{Mod}_{\Lambda}$ to the category of $\Lambda$-modules.

%

\begin{prop}\label{comp} We have a bifunctorial isomorphism
$$
\Gamma(U,T(\F)) \xrightarrow[\sim]{} T(\Gamma(U,\F)),
$$
for $U$ in $\Et_{Y_{\bs}}$ and $\F$ in $\Sh(Y,Q, \Lambda)$.
\end{prop}

Ler us consider a factorization $\overline{s} \xrightarrow[]{\varphi_0} s' \rightarrow s$ where $s'$ is a finite extension of $s$, and an object $U'$ of $\Et_{Y_{s'}}$, such that $U$ is isomorphic to $U' \times_{s', \varphi_0} \bs$. We henceforth assume that $U$ is equal to $U' \times_{s',\varphi_0} \bs$. We have a decomposition
$$
U' \times_{s} \bs =  \coprod_{\overline{s} \xrightarrow[]{\varphi} s'} U' \times_{s',\varphi} \bs,
$$
hence a splitting
$$
\F(U'_{\bs}) \cong \bigoplus_{\overline{s} \xrightarrow[]{\varphi} s'} \F(  U' \times_{s',\varphi} \bs),
$$
and in particular
$$
T\left( \F(U'_{\bs}) \right) \cong \bigoplus_{\overline{s} \xrightarrow[]{\varphi} s'} T \left( \F(  U' \times_{s',\varphi} \bs) \right).
$$
Taking into consideration the action of $G_s$ on the left hand side, cf. \ref{5.6}, we obtain an isomorphism
$$
T_{G_s}\left( \F(U'_{\bs}) \right) \cong \Ind_{G_{s'}}^{G_s} T_{G_{s'}} \left( \F(  U' \times_{s',\varphi_0} \bs) \right).
$$
Taking the $G_s$-invariants, this yields
$$
\Gamma(U',T(\F)) \cong  T_{G_{s'}} \left( \F(  U' \times_{s',\varphi_0} \bs) \right)^{G_{s'}} = T_{G_{s'}} \left( \F( U) \right)^{G_{s'}}
$$
If $\overline{s} \rightarrow s'' \rightarrow s'$ is a factorization of $\varphi_0$, where $s''$ is a finite extension of $s'$, then replacing $s'$ and $U'$ by $s''$ and $U' \times_{s'} s''$ above yields
$$
\Gamma(U'\times_{s'} s'',T(\F)) \cong   T_{G_{s'}} \left( \F( U) \right)^{G_{s''}},
$$
By taking the colimit over $s''$, we obtain
$$
\Gamma(U,T(\F)) \cong   T \left( \F( U) \right).
$$

\begin{cor}\label{isexact} For any geometric point $\bar{y}$ of $Y$, we have a natural isomorphism
$$
T(\F)_{\bar{y}} \cong T( \F_{\bar{y}}),
$$
for $\F$ in $\Sh(Y,Q, \Lambda)$. In particular, the functor 
$$
T : \Sh(Y,Q, \Lambda) \rightarrow \Sh(Y, \Lambda),
$$
is exact and commutes with small filtered colimits.
\end{cor}

Indeed, we can consider $\bar{y}$ as a geometric point of $Y_{\bs}$ and the conclusion follows by taking the (filtered) colimit of the isomorphisms
$$
\Gamma(U,T(\F)) \xrightarrow[\sim]{} T(\Gamma(U,\F)),
$$
from \ref{comp}, where $U$ runs over qcqs \'etale neighbourhoods of $\bar{y}$ over $Y_{\bs}$.

\begin{cor}\label{pullback} Let $f : Y' \rightarrow Y$ be a morphism of $s$-schemes. For any object $\F$ of $\Sh(Y,Q, \Lambda)$, the natural homomorphism
$$
f^{-1} T(\F) \rightarrow T(f^{-1} \F),
$$
with notation as in \ref{5.2}, is an isomorphism.
\end{cor}

This is an immediate consequence of \ref{isexact}.


\begin{prop}\label{acyclicity2} Let $\F$ be an injective object in $\Sh(Y,Q, \Lambda)$. Then for any $U$ in $\Et_{Y_{\bs}}$ and any $j > 0$, we have $H^j(U,T(F)) = 0$. 
\end{prop}

For any \'etale cover $\mathcal{U} = (U_i \rightarrow U)_{i \in I}$ in $\Et_{Y_{\bs}}$ with $I$ finite, the \v{C}ech complex $\check{\mathcal{C}}^\bullet (\mathcal{U}, \F)$ is acyclic in nonzero degree by Proposition \ref{acyclicity1}. By Proposition \ref{comp} and by exactness of $T$, we obtain that
$$
\check{\mathcal{C}}^\bullet (\mathcal{U}, T(\F)) \cong T\left(\check{\mathcal{C}}^\bullet (\mathcal{U}, \F) \right),
$$
is acyclic as well in nonzero degree. The conclusion then follows from (\cite{stacks-project} 03F9).


\subsection{\label{5.4}} Let $\theta,\Lambda$ and $(T,\tau)$ be as in \ref{5.7}. Let $f : Y' \rightarrow Y$ be a qcqs morphism of $s$-schemes. For any object $U$ of $\Et_Y$, the fiber product $U \times_{Y} Y'$ is an object of $\Et_{Y'}$. Moreover, for any object $\F$ of $\Sh(Y,Q, \Lambda)$, we have, with the notation from \ref{5.2},
\begin{align*}
\Gamma(U,T(f_*\F)) &= T_{G_s} \left( f_*\F(U_{\bs} \right)^{G_s} \\
&=  T_{G_s} \left( \F(U_{\bs} \times_{Y_{\bs}} Y'_{\bs} \right)^{G_s} \\
&=  \Gamma(U,f_*T(\F)),
\end{align*}
hence a commutative diagram

\begin{center}
 \begin{tikzpicture}[scale=1]

\node (A) at (0,0) {$\Sh(Y,Q, \Lambda)$};
\node (B) at (3,0) {$\Sh(Y, \Lambda)$};
\node (C) at (0,2) {$\Sh(Y',Q, \Lambda)$};
\node (D) at (3,2) {$\Sh(Y', \Lambda)$};

\path[->,font=\scriptsize]
(A) edge node[above]{$T$} (B)
(D) edge node[left]{$f_*$}  (B)
(C) edge node[left]{$f_*$} (A)
(C) edge node[above]{$T$} (D);
\end{tikzpicture} 
\end{center}
of abelian categories.

\begin{prop}\label{pushforwardderived} We have a natural isomorphism
$$
Rf_* \circ T \cong T \circ Rf_*,
$$
of functors from $D^+(Y',Q, \Lambda)$ to $D^+(Y,\Lambda)$.
\end{prop}

Indeed, for any injective object $\F$ of $\Sh(Y',Q, \Lambda)$, the object $f_* \F$ of $\Sh(Y,Q, \Lambda)$ is $T$-acyclic by exactness of $T$, cf. \ref{isexact}. Moreover, for any integer $j > 0$ the restriction of $R^j f_* T(\F)$ to $Y_{\bs}$ is the sheaf associated to the presheaf
$$
U \mapsto H^j(U \times_{Y_{\bs}} Y_{\bs}', T(\F)),
$$
which vanishes by \ref{acyclicity2}. Thus $T(\F)$ is $f_*$-acyclic. This implies that the canonical transformations
\begin{align*}
R(f_* \circ T) &\rightarrow Rf_* \circ T \\
R( T \circ f_*) &\rightarrow  T \circ Rf_*,
\end{align*}
are isomorphisms, and the conclusion then follows from the equality $f_* \circ T = T \circ f_*$. 

\subsection{\label{5.9}} Let $\theta : Q \rightarrow G_s$ be a surjective continuous homomorphism of profinite groups with kernel $I$ and let $\Lambda$ be a noetherian ring. Let $Y$ be a locally noetherian $s$-scheme. We denote by $\Sh_c(Y,Q,\Lambda)$ the full subcategory of $\Sh(Y,Q,\Lambda)$ consisting of objects $\F$ whose underlying sheaf of $\Lambda$-module on $Y_{\bs}$ is constructible, cf. (\cite{SGA4}, IX 2.3).

\begin{prop}\label{constr10} Let $(T,\tau)$ be a $\Lambda$-linear $(Q,I)$-contraction and let $\F$ be an object of $\Sh(Y,Q,\Lambda)$.
\begin{enumerate}
\item If the object of $\Sh(Y_{\bs},\Lambda)$ underlying $\F$ is locally constant constructible, then $T(\F)$ is locally constant constructible.
\item If the object of $\Sh(Y_{\bs},\Lambda)$ underlying $\F$ is constructible, then $T(\F)$ is constructible.
\end{enumerate}
\end{prop}

Part (1) follows from Corollary \ref{isexact}, from \ref{contract} (2), and from (\cite{SGA4}, IX 2.13 (i)). Part (2) follows from (1) and from Corollary \ref{pullback}.

\section{Nearby and vanishing cycles \label{7}}

Let $S$ be a scheme and let $s$ be a point of $S$ such that the local ring $\Ow_{S,s}$ is a discrete valuation ring. We denote by $\bs$ a separable closure of $s$, by $\eta$ the generic point of the henselization $S_{(s)}$ of $S$ at $s$, and by $\bta$ a separable closure of the generic point $\eta_{\bs}$ of the strict henselization $S_{(\bs)}$.

\subsection{\label{6.1}} Let us denote by $G_{\eta}$ and $G_s$ the Galois group of the extensions $\bta \rightarrow \eta$ and $\bs \rightarrow s$ respectively. We have a natural specialization homomorphism $\theta : G_{\eta} \rightarrow G_s$. This homomorphism is surjective and we denote by $I_{\eta}$ its kernel, so that we have an exact sequence
$$
1 \rightarrow I_{\eta} \rightarrow G_{\eta} \rightarrow G_s \rightarrow 1,
$$
of profinite groups. In this section, we apply the definitions and results of Section \ref{5} to $Q = G_{\eta}$ and $\theta$ as above.

\subsection{\label{6.2}} Let $f : X \rightarrow S$ be an $S$-scheme. Let us consider the commutative diagram

\begin{center}
 \begin{tikzpicture}[scale=2]

\node (E) at (0,0) {$\bs$};
\node (F) at (1,0) {$ S_{(\bs)}$};
\node (G) at (2,0) {$ S_{(s)}$};
\node (H) at (3,0) {$S$};
\node (A) at (0,1) {$X_{\bs}$};
\node (B) at (1,1) {$X \times_{S} S_{(\bs)}$};
\node (C) at (2,1) {$X \times_{S} S_{s}$};
\node (D) at (3,1) {$X$};

\path[->,font=\scriptsize]
(E) edge  (F)
(F) edge  (G)
(G) edge  (H)
(A) edge node[above]{$\overline{i}$} (B)
(B) edge  (C)
(C) edge   (D)
(A) edge node[ right]{$f_{\bs}$} (E)
(B) edge  (F)
(C) edge  (G)
(D) edge node[ left]{$f$} (H);
\end{tikzpicture} 
\end{center}
of $S$-schemes, where all squares are cartesian.

\begin{defi}[\cite{SGA7}, XIII 1.3.2] Let $\Lambda$ be a ring and let $\F$ be an object of $\Sh(X,\Lambda)$. We denote by $\Psi_f^s(\F)$ the object of $\Sh(X_s,G_{\eta},\Lambda)$, cf. \ref{equivsheaf}, obtained the pullback by $\overline{i}$ of the sheaf of $\Lambda$-modules on $ X \times_{S} S_{(\bs)}$ whose sections on an object $U$ of $\Et_{ X \times_{S} S_{(\bs)}}$ are given by $\Gamma(U \times_{S_{(\bs)}} \bta, \F)$.
\end{defi}

This defines a left exact functor
$$
\Psi_f^s : \Sh(X,\Lambda) \rightarrow \Sh(X_s,G_{\eta},\Lambda),
$$
thereby inducing a right derived functor 
$$
R \Psi_f^s : D^+(X,\Lambda) \rightarrow D^+(X_s,G_{\eta},\Lambda),
$$
namely the \textit{nearby cycles} functor, cf. (\cite{SGA7}, XIII 2.1.1). 

\subsection{\label{6.3}} Let $i : X_s \rightarrow X$ be the natural closed immersion. The natural cospecialization homomorphism
$$
i^{-1} \rightarrow \Psi_f^s,
$$
induces a natural transformation
$$
i^{-1} \rightarrow R\Psi_f^s,
$$
of triangulated functors from $D^+(X,\Lambda)$ to $D^+(X_s,G_{\eta},\Lambda)$. We denote its cone by
$$
R \Phi_f^s  : D^+(X,\Lambda) \rightarrow D^+(X_s,G_{\eta},\Lambda),
$$
the \textit{vanishing cycles} functor. We thus have a functorial distinguished triangle
$$
i^{-1} \F \rightarrow R\Psi_f^s \F \rightarrow R \Phi_f^s(\F) \xrightarrow[]{[1]},
$$
for $\F$ in $D^+(X,\Lambda)$.

\subsection{\label{6.4}} Let $\Lambda$ be a noetherian ring annihilated by an integer invertible on $s$. Let $(T,\tau)$ be a $\Lambda$-linear $(G_{\eta},I_{\eta})$-contraction, cf. \ref{contract}. The composition of the functor $T$ from \ref{5.8} with the vanishing cycles functor from \ref{6.3} yields a triangulated functor
$$
T R \Phi_f^s  : D^+(X,\Lambda) \rightarrow D^+(X_s,\Lambda).
$$
We now reference a few of its properties.

\begin{prop}\label{constr11} Let $\F$ be an object of $D_c^b(X,\Lambda)$. If $f$ is locally of finite type, then $T R \Phi_f^s$ belongs to $D_c^b(X_s,\Lambda)$.
\end{prop}

Under these assumptions, it follows from (\cite{SGA412}, 7.3.2) that $R\Psi_{f}^s(\F)$ belongs to $D_c^b(X_s,\Lambda)$, hence $R\Phi_{f}^s(\F)$ belongs as well to $D_c^b(X_s,\Lambda)$ and so does $T R \Phi_f^s(\F)$ by \ref{constr10} (2).

\begin{prop}\label{isom5} Let $h : X' \rightarrow X$ be a smooth morphism of $S$-schemes and let $\F$ be an object of $D^+(X,\Lambda)$. Then the natural homomorphism
$$
h_s^{-1}  T R \Phi_f^s(\F)  \rightarrow T R \Phi_{fh}^s(h^{-1}\F),
$$
is an isomorphism, where $h_s : X'_s \rightarrow X_s$ is the morphism induced by $h$. In particular, if $f$ is smooth then $T R \Phi_f^s(\Lambda)$ vanishes.
\end{prop}

This follows from the corresponding property of $R\Psi_{f}^s$, cf. (\cite{SGA7}, XIII 2.1.7.1), and from Corollary \ref{pullback}.

\begin{prop}\label{isom6} Let $h : X' \rightarrow X$ be a proper morphism of $S$-schemes and let $\F$ be an object of $D^+(X',\Lambda)$. Then we have a natural homomorphism
$$
Rh_{s*}  T R \Phi_{fh}^s(\F')  \rightarrow T R \Phi_f^s( Rh_* \F),
$$
in $D^+(X_s,\Lambda)$. In particular, if $f$ is proper then we have a functorial isomorphism
$$
R\Gamma(X_{\overline{s}}, T R \Phi_f^s(\F)) \rightarrow  T R \Phi_{\id}^s( Rf_* \F)
$$
in $D^+(G_s-\Mod_{\Lambda})$.
\end{prop}

This follows from the corresponding property of $R\Psi_{f}^s$, cf. (\cite{SGA7}, XIII 2.1.7.2), and from Proposition \ref{pushforwardderived} since $h$ is qcqs.

\subsection{\label{6.5}} For any finite local $\mathbb{Z}_{\ell}$-algebra, let $(T_{\Lambda},\tau_{\Lambda})$ be a $\Lambda$-linear $(G_{\eta},I_{\eta})$-contraction, cf. \ref{contract}, together with compatible isomorphisms
$$
T_{\Lambda} \otimes_{\Lambda} \Lambda' \cong T_{\Lambda'},
$$
for any homomorphism $\Lambda \rightarrow \Lambda'$ of finite local $\mathbb{Z}_{\ell}$-algebra.

\begin{defi}\label{extensionadm} Let $\Lambda$ be an admissible $\ell$-adic coefficient ring, cf. (\cite{G19}, 2.2). Let $M$ be a finitely generated $\Lambda$-module endowed with a continuous $\Lambda$-linear left $I_{\eta_s}$-action.
\begin{enumerate}
\item If $\Lambda$ is the ring of integers in a finite extension of $\mathbb{Q}_{\ell}$, with maximal ideal $\m$, we set
$$
T_{\Lambda}(M) = \lim_{n} T_{\Lambda / \m^n}(M/\m^n). 
$$
\item If $\Lambda$ is a finite extension of $\mathbb{Q}_{\ell}$, with ring of integers $\Lambda_0$, and if $M_0$ is an $I_{\eta_s}$-stable finitely generated sub-$\Lambda_0$-module such that $\Lambda M_0 = M$, we set
$$
T_{\Lambda}(M) =  T_{\Lambda_0}(M_0) \otimes_{\Lambda_0} \Lambda.
$$
This is independent of the choice of $M_0$.
\end{enumerate}
\end{defi}

Let $f : X \rightarrow S$ be a morphism locally of finite type, with $\ell$-invertible on $S$. Proposition \ref{constr11} allows to extend the functors from \ref{extensionadm} to triangulated functors
$$
T R \Phi_f^s  : D_c^b(X,\Lambda) \rightarrow D_c^b(X_s,\Lambda),
$$
for any admissible $\ell$-adic coefficient ring, in such a way that its formation commutes with base change along homomorphisms $\Lambda \rightarrow \Lambda'$ of admissible $\ell$-adic coefficient rings, and so that it coincides with the functor $T R \Phi_f^s$ from \ref{6.4} when $\Lambda$ is a finite $\mathbb{Z}_{\ell}$-algebra. Propositions \ref{isom5} and \ref{isom6} still hold in this context.

\section{Gabber-Katz extensions and examples of contractions \label{0}}

In this section, we provide examples of linear contractions, cf. \ref{contract}, whose construction are based on the theory of Gabber-Katz extensions from \cite{katz}. We review the latter in \ref{3.1} by placing these extensions in what we believe to be their rightful context, namely as pullbacks along a morphism of toposes, the \textit{canonical retraction}, cf \ref{3.1.3}.  Throughout this section, we fix a perfect field $k$ of characteristic $p > 0$.

\subsection{\label{3.1}} Let $S$ be a smooth $k$-scheme purely of dimension $1$, and let $s$ be a closed point of $S$. Thus the henselization $S_{(s)}$ of $S$ at $s$ is a henselian trait whose residue field $k(s)$ is a finite extension of $k$. Let $\pi$ be a uniformizer of $S_{(s)}$, thereby inducing a morphism
$$
\pi : S_{(s)} \rightarrow \mathbb{A}^1_s,
$$
corresponding to the unique homomorphism $k(s)[t] \rightarrow \Ow_{S_{(s)}}$ of $k(s)$-algebras sending $t$ to $\pi$. Let us denote by 
$$
\pi_{\diamondsuit} : \Et_{S_{(s)}} \rightarrow \Et_{\mathbb{A}^1_s},
$$
the Gabber-Katz extension functor, cf. (\cite{G19}, 4.2): for any \'etale $S_{(s)}$-scheme $U$, the \'etale $\mathbb{A}^1_s$-scheme $\pi_{\diamondsuit} U$ is the unique (up to isomorphism) such $\mathbb{A}^1_s$-scheme whose pullback by $\pi$ is isomorphic to $U$ and whose restriction along the open immersion $\mathbb{G}_{m,s} \rightarrow \mathbb{A}^1_s$ is a \textit{special} finite \'etale $\mathbb{G}_{m,s}$-scheme, namely it is tamely ramified above $\infty$ and its geometric monodromy group has a unique $p$-Sylow, cf. (\cite{G19}, 4.3). We will also label \textit{special} an \'etale $\mathbb{A}^1_s$-scheme whose restriction along the open immersion $\mathbb{G}_{m,s} \rightarrow \mathbb{A}^1_s$ is \textit{special}, or equivalently an object of the essential image of $\pi_{\diamondsuit}$.

\begin{prop}\label{3.1.1} The functor $\pi_{\diamondsuit}$ has the following properties:
\begin{enumerate}
\item it sends a final object of $\Et_{S_{(s)}}$ to a final object of $\Et_{\mathbb{A}^1_s}$;
\item for any covering $(U_{i } \rightarrow U)_{i \in I}$ of the site $\Et_{S_{(s)}}$, the family $(\pi_{\diamondsuit} U_{i} \rightarrow \pi_{\diamondsuit} U)_{i \in I}$ is a covering in site $\Et_{\mathbb{A}^1_s}$;
\item for any diagram $U \rightarrow W \leftarrow V$ in $\Et_{S_{(s)}}$, the natural $\mathbb{A}^1_s$-morphism
$$
\pi_{\diamondsuit}(U \times_W V) \rightarrow \pi_{\diamondsuit}(U) \times_{\pi_{\diamondsuit}(W)} \pi_{\diamondsuit}(V),
$$
is an isomorphism.
\end{enumerate}
\end{prop}

The $\mathbb{A}^1_s$-scheme $\mathbb{A}^1_s$ is special and its pullback by $\pi$ is $S_{(s)}$, hence $\pi_{\diamondsuit} S_{(s)}$ isomorphic to $\mathbb{A}^1_s$. This proves $(1)$. In order to prove $(2)$, we can assume that $U$ is connected, since the functor $\pi_{\diamondsuit}$ commutes with finite disjoint unions. Moreover, we have that $\pi^{-1} \pi_{\diamondsuit} U_{i}$ is isomorphic to $ U_{i}$, so that $(\pi^{-1} \pi_{\diamondsuit} U_{i})_{i \in I}$ covers $\pi^{-1} \pi_{\diamondsuit} U$, hence is is sufficient to prove that $(j^{-1} \pi_{\diamondsuit} U_{i})_{i \in I}$ covers $j^{-1} \pi_{\diamondsuit} U$, where $j : \mathbb{G}_{m,s} \rightarrow \mathbb{A}^1_s$ is the natural open immersion. Let $i_0$ in $I$ be such that the image of $U_{i_0}$ in $U$ contains the generic point of $U$. Then the image of $j^{-1}\pi_{\diamondsuit} U_{i_0}$ in $j^{-1}\pi_{\diamondsuit} U$ is non empty. Since the morphism $j^{-1}\pi_{\diamondsuit} U_{i_0} \rightarrow j^{-1}\pi_{\diamondsuit} U$ is finite \'etale, its image is closed and open in the connected scheme $j^{-1}\pi_{\diamondsuit} U$. Consequently, the morphism $j^{-1}\pi_{\diamondsuit} U_{i_0} \rightarrow j^{-1}\pi_{\diamondsuit} U$ is surjective.
Let us now prove $(3)$. Let 
$$
\alpha : 
\pi_{\diamondsuit}(U \times_W V) \rightarrow \pi_{\diamondsuit}(U) \times_{\pi_{\diamondsuit}(W)} \pi_{\diamondsuit}(V),
$$
be the morphism under consideration. Then $\pi^{-1} \alpha$ is an isomorphism, and we infer from the characterization of $\pi_{\diamondsuit}$, cf. \ref{3.1}, that $(3)$ holds if and only if $X = \pi_{\diamondsuit}(U) \times_{\pi_{\diamondsuit}(W)} \pi_{\diamondsuit}(V)$ is a special $\mathbb{A}^1_s$-scheme. Since $\pi_{\diamondsuit}(U) $ and $\pi_{\diamondsuit}(W)$ are tamely ramified above infinity, so is $X$. Moreover, the geometric monodromy group of $j^{-1} X \rightarrow \mathbb{G}_{m,s}$ is a subgroup of the product of those of $j^{-1}\pi_{\diamondsuit}(U) $ and $j^{-1}\pi_{\diamondsuit}(U) $, and the class of finite groups which admit a unique $p$-Sylow is stable by finite products and subobjects, hence the geometric monodromy group of $j^{-1} X \rightarrow \mathbb{G}_{m,s}$ belongs to that class. Thus $X$ is a special $\mathbb{A}^1_s$-scheme, and $(3)$ holds.

\begin{cor}\label{3.1.2} The functor $\pi_{\diamondsuit} : \Et_{S_{(s)}} \rightarrow \Et_{\mathbb{A}^1_s}$ is a morphism of sites.
\end{cor}

This follows from Proposition \ref{3.1.1} and from the fact that $\Et_{S_{(s)}}$ is an essentially small category.

\begin{defi}\label{3.1.3}Let $(S,s,\pi)$ be as in \ref{3.1}. The \textit{Gabber-Katz retraction} associated to this data is the morphism of toposes $r_{s,\pi} : \mathbb{A}^1_{s,\et} \rightarrow S_{(s),\et}$ associated to the morphism of sites $\pi_{\diamondsuit}$, cf. \ref{3.1.2}.
\end{defi}

The composition of $r_{s,\pi}$ with the morphism $\pi : S_{(s),\et} \rightarrow \mathbb{A}^1_{s,\et} $ of toposes induced by $\pi$, is (equivalent to) the identity morphism of $S_{(s),\et}$, hence the terminology ``retraction''. The pair $r_{s,\pi} = (r_{s,\pi}^{-1}, r_{s,\pi*})$ of adjoint functors can be described as follows. For an object $\F$ of $\mathbb{A}^1_{s,\et}$, the object $r_{s,\pi*} \F$ is the sheaf on $\Et_{S_{(s)}}$ given by $U \mapsto \F( \pi_{\diamondsuit} U)$. Moreover, the functor $r_{s,\pi}^{-1} : S_{(s),\et} \rightarrow \mathbb{A}^1_{s,\et}$ is the Gabber-Katz extension functor, sending a representable sheaf $\underline{U}$ to the representable sheaf $\underline{\pi_{\diamondsuit} U}$.

\subsection{\label{3.2}} Let $\Lambda$ be an admissible $\ell$-adic coefficient ring, cf. (\cite{G19}, 2.2), and let $(S,s,\pi)$ be as in \ref{3.1}. Let $\eta_s$ be the generic point of $S_{(s)}$, let $\bs$ be a separable closure of $s$, and let $\bta_s$ be a separable closure of $\eta_{\bs}$. Let $G_{\eta_s}$ and $G_s$ be the Galois group of the extensions $\bs \rightarrow s$ and $\bta_{s} \rightarrow \eta_s$. Let $I_{\eta_s}$ be the kernel of the sujective specialization homomorphism $G_{\eta_s} \rightarrow G_s$. For an object $M$ of $G_{\eta_s}-\Mod_{\Lambda}$, we denote by $M_!$ the extension by zero to $S_{(s)}$ of the corresponding sheaf on $\eta_s$. In (\cite{katz}, 1.6), Katz considered the functor defined by 
$$
\mathrm{Sw}(M) = H^1( \mathbb{A}^1_{\bs},  r_{\bs,\pi}^{-1} M_!),
$$
for $M$ in $I_{\eta_s}-\Mod_{\Lambda}$, where $r_{\bs,\pi}$ is the canonical retraction from \ref{3.1.3}. This is naturally a $\Lambda$-linear $(G_{\eta_s},I_{\eta_s})$-retraction if $\Lambda$ is finite, and yields a compatible system as in \ref{6.5} in general. When $\Lambda$ is a field, this functor has the remarkable property that 
$$
\rk(\mathrm{Sw}(M)) = \mathrm{sw}(M),
$$
where $\mathrm{sw}(M)$ is the Swan conductor of $M$, hence $\mathrm{Sw}$ can be considered as a \textit{refined Swan conductor}. This particular $\Lambda$-linear contraction will not be appear in the rest of this text.

\subsection{\label{3.3}} We retain the notation from \ref{3.2}. Let us set $\tilde{\pi} = (1-\pi^{-1})^{-1}$. The uniformizer $\pi$ is uniquely determined by $\tilde{\pi}$, since we have $ \pi = (1-\tilde{\pi}^{-1})^{-1}$. Let $\overline{1} : \bs \rightarrow \mathbb{A}^1_{\bs}$ be the $\bs$-point $1$ of $\mathbb{A}^1_{\bs}$. Let us set
$$
\mathrm{CFT}_{\tilde{\pi}}(M) = r_{\bs, \pi}^{-1}(M_!)_{\overline{1}},
$$
for $M$ in $I_{\eta_s}-\Mod_{\Lambda}$. This is naturally a $\Lambda$-linear $(G_{\eta_s},I_{\eta_s})$-retraction if $\Lambda$ is finite, and yields a compatible system as in \ref{6.5} in general. The functor $\mathrm{CFT}_{\tilde{\pi}}$ can be considered as constituting a \textit{refined class field theory} for $\eta_s$, as indicated by the following result.

\begin{prop}\label{CFT} Let us assume that $\Lambda$ is a field. Let $M$ be an object of $G_{\eta_s}-\Mod_{\Lambda}$ and let $D$ be a divisor on $S_{(s)}$ such that $\det(M)$ has ramification bounded $D$, cf. $($\cite{G19}, $5.17)$. Let $\chi_{\det(M)}$ be the multiplicative local system on $\Pic(S_{(s)},D)_s$, cf. $($\cite{G19}, $5.20)$, associated to $\det(M)$ by geometric class field theory, cf. $($\cite{G19}, $5.25)$. Then we have
$$
\langle \chi_{\det(M)} \rangle \left( \tilde{\pi} \right) = \det(\mathrm{CFT}_{\tilde{\pi}}(M) ),
$$
with the notation from $($\cite{G19}, $7.11)$.
\end{prop}

Indeed, we have
$$
\det( \mathrm{CFT}_{\tilde{\pi}}(M) ) = \mathrm{CFT}_{\tilde{\pi}}(\det(M)),
$$
and (\cite{G19}, 5.34) yields
$$
(\chi_{\det(M)})_{| 1 -  \pi^{-1}} = \det( \mathrm{CFT}_{\tilde{\pi}}(M) ),
$$
hence the conclusion.

\subsection{\label{0.1}} We still retain the notation from \ref{3.2}. Let $\Lambda$ be an admissible $\ell$-adic coefficient ring, cf. (\cite{G19}, 2.2). The profinite group $I_{\eta_s}$ admits a unique $p$-Sylow subgroup $P$, the \emph{wild inertia group} associated to $(S_{(s)},\overline{\eta_s})$.

\begin{defi}\label{totram} A continuous representation of $G_{\eta_s}$ on a $\Lambda$-module $M$ is \emph{totally wildly ramified} if the submodule $M^P$ of $P$-invariants elements vanishes. 
\end{defi}

\begin{exemple}\label{totramex} Let $\psi : \mathbb{F}_p \rightarrow \Lambda^{\times}$ be a non trivial homomorphism and let $f$ be an element of $k(\eta_s)$ whose valuation is negative and prime to $p$. Then the continuous representation of $G_{\eta}$ corresponding to $\Lc_{\psi} \lbrace f \rbrace$ is totally wildly ramified. The latter assertion holds as well for $\F \otimes \Lc_{\psi} \lbrace f \rbrace$ for any tamely ramified lisse \'etale sheaf of $\Lambda$-modules $\F$ on $\eta$, since $P$ acts trivially on $\F_{\overline{\eta}}$.
\end{exemple}

\begin{prop}\label{vanish1} Let $M$ be a $\Lambda$-module endowed with a continuous $\Lambda$-linear left action of $G_{\eta}$, and let $\F$ be the corresponding lisse \'etale sheaf of $\Lambda$-modules on $\eta$. If $M$ is totally wildly ramified (cf. \ref{totram}) then the canonical homomorphism
$$
j_! \F \rightarrow R j_* \F,
$$ 
is a quasi-isomorphism.
\end{prop}

This amounts to the acyclicity of the fiber of $R j_* \F$ at $\overline{s}$, i.e. to the vanishing of $H^j(I_{\eta_s}, M)$ for each integer $j$. We have the Hochschild-Serre spectral sequence
$$
H^i( I_{\eta_s}/P, H^j(P, M)) \Rightarrow H^{i+j}(I_{\eta_s},M).
$$
The groups $H^j(P, M)$ vanish for $j \neq 0$ since $P$ is a $p$-group and $M$ is of $\ell$-torsion, while $H^0(P,M) = M^P$ vanishes since $M$ is assumed to be totally wildly ramified. Thus $H^j(P, M)$ vanishes for any $j$, and consequently so does $H^j(I_{\eta_s}, M)$.

\begin{prop}\label{isom1} Let $M$ be a $\Lambda$-module endowed with a continuous $\Lambda$-linear left action of $I_{\eta_s}$. Then the canonical homomorphism
$$
R \Gamma( \mathbb{A}^1_{\bs}, r_{\bs,\pi}^{-1} M_! \otimes \Lc_{\psi} \lbrace -t \rbrace) \rightarrow R \Gamma_c( \mathbb{A}^1_{\bs}, r_{\bs,\pi}^{-1} M_! \otimes \Lc_{\psi} \lbrace -t \rbrace),
$$
is a quasi-isomorphism. Moreover, the complex $R \Gamma( \mathbb{A}^1_{\bs},r_{\bs,\pi}^{-1} M_! \otimes \Lc_{\psi} \lbrace -t \rbrace)$ is concentrated in degree $1$.
\end{prop}

The first assertion follows from Example \ref{totramex} and from Proposition \ref{vanish1}. The second assertion follows from the vanishing of $H^0(\mathbb{A}^1_{\bs},r_{\pi}^{-1} M_! \otimes \Lc_{\psi} \lbrace -t \rbrace)$ and from the fact that the \'etale cohomological dimension of a smooth affine curve over a separably closed field is at most $1$. 

\begin{cor}\label{exact1} The functor
$$
\Art_{\pi} : M \mapsto H_c^1(  \mathbb{A}^1_{\bs}, r_{\bs,\pi}^{-1} M_!  \otimes \Lc_{\psi} \lbrace -t \rbrace),
$$
from $I_{\eta_s}-\Mod_{\Lambda}$ to $\Mod_{\Lambda}$, is additive, exact and commutes with base change along homomorphisms $\Lambda \rightarrow \Lambda'$ of admissible $\ell$-adic coefficient rings.
\end{cor}

This is an immediate consequence of Proposition \ref{isom1}. We note that if $\Lambda$ is finite, then $\Art_{\pi}$ is naturally endowed with a structure of $\Lambda$-linear $(G_{\eta_s},I_{\eta_s})$-contraction, and we obtain a compatible system as in \ref{6.5}.

%
\begin{prop}\label{art4} Let assume that $\Lambda$ is a field. Let $M$ be a $\Lambda$-module endowed with a continuous $\Lambda$-linear left action of $G_{\eta_s}$. Then we have
\begin{align*}
\rk( \Art_{\pi}(M) )&= \rk(M) + \mathrm{sw}(M) = a(S_{(s)},M_!), \\
\det( \Art_{\pi}(M) )&= \varepsilon(S_{(s)},M_!,d\pi),
\end{align*}
where $a(S_{(s)},-)$ is the Artin conductor, cf. (\cite{G19}, 7.2), and $\varepsilon(S_{(s)},-)$ is the geometric local $\varepsilon$-factor from (\cite{G19}, 9.2).
\end{prop}

The first assertion follows from the Grothendieck-Ogg-Shafarevich formula, while the second one matches the definition of geometric local $\varepsilon$-factors in (\cite{G19}, 9.2).

\begin{cor}\label{art10} Let assume that $\Lambda$ is a field, and let $\F$ be an object of $D_c^b(S_{(s)},\Lambda)$. Then we have 
\begin{align*}
\rk( \Art_{\pi} R\Phi_{\id}^s(\F) )&= a(S_{(s)},\F), \\
\det( \Art_{\pi} R\Phi_{\id}^s(\F) )&= \varepsilon(S_{(s)},\F,d\pi).
\end{align*}
\end{cor}

This follows from Propositions \ref{art4} and \ref{tensor}, and from the fact that $\Art_{\pi}(\Lambda) \cong \Lambda$. 

\section{Proof of Theorem \ref{teo0} \label{prf} }

In this section, we prove Theorem \ref{teo0}. The preliminary paragraph \ref{4.0} provides a more detailed conclusion in the case of a projective line. We then explain how to reduce to projective schemes by using Chow's lemma in \ref{4.1}. Theorem \ref{teo0} is then reduced to the case of projective spaces in \ref{4.2}, and the latter case is handled by induction on the dimension in \ref{4.3}. Throughout this section, we denote by $\Lambda$ a finite extension of either $\mathbb{F}_{\ell}$ of $\mathbb{Q}_{\ell}$.

\subsection{\label{4.0}} Let $S = \mathbb{P}^1_k$, and let $\F$ be an object of $D_c^b(S, \Lambda)$. We denote by $t$ the canonical coordinate on the open subset $\mathbb{A}^1_k$ of $S$. By the product formula (\cite{G19}, 1.4), we have,  with notation as in (\cite{G19}, 3.26),
\begin{align*}
\varepsilon_{\bk}(S, \F) &= \chi_{\cyc}^{-\rk(\F)} \varepsilon(S_{(\infty)},\F_{|S_{(\infty)}},dt) \prod_{s \in |\mathbb{A}^1_k|} \delta_{s/k}^{a(S_{(s)},\F_{|S_{(s)}})} \Ver_{s/k} \left( \varepsilon(S_{(s)},\F_{|S_{(s)}},dt) \right) \\
&=  \chi_{\cyc}^{\rk(\F)} \langle \chi_{\det(\F_{\eta_{\infty}})} \rangle(-t^2) \varepsilon(S_{(\infty)},\F_{|S_{(\infty)}},d(t^{-1})) \prod_{s \in |\mathbb{A}^1_k|} \delta_{s/k}^{a(S_{(s)},\F_{|S_{(s)}})} \Ver_{s/k} \left( \varepsilon(S_{(s)},\F_{|S_{(s)}},dt) \right) 
\end{align*}
For any closed point $s$ of $\mathbb{A}^1_k$, we have $dt = d(t-s)$ and $\pi_s = t-s$ is a uniformizer of $S_{(s)}$. For $s = \infty$, we set $\pi_{\infty} = t^{-1}$. Then Proposition \ref{art10} yields
$$
 \delta_{s/k}^{a(S_{(s)},\F_{|S_{(s)}})} \Ver_{s/k} \left( \varepsilon(S_{(s)},\F_{|S_{(s)}},d\pi_s) \right)  = \det \left( \Ind_{G_s}^{G_k} \Art_{\pi_s} R\Phi_{\id}^s(\F) \right),
$$
for any  closed point $s$ of $\mathbb{P}^1_k$. By Proposition \ref{CFT}, we further have
$$
 \langle \chi_{\det(\F_{\eta_{\infty}})} \rangle(-t^2) = \det \left( \mathrm{CFT}_{t^{-1}}(\F_{\eta_{\infty}}) \right)^{-1}  \det \left( \mathrm{CFT}_{-t^{-1}}(\F_{\eta_{\infty}}) \right)^{-1}.
$$
We thus obtain that $\varepsilon_{\bk}(S, \F)$ coincides with
\begin{align*}
 \det \left( \mathrm{CFT}_{t^{-1}} R \Psi_{\id}^{\infty} \F(-1) \right)^{-1}  \det \left( \mathrm{CFT}_{-t^{-1}}  R \Psi_{\id}^{\infty} \F(-1) \right)^{-1} \prod_{s \in |\mathbb{P}^1_k|} \det \left( \Ind_{G_s}^{G_k} \Art_{\pi_s} R\Phi_{\id}^s(\F) \right) .
\end{align*}
%

\subsection{\label{4.1}} Let us assume that Theorem \ref{teo0} is proved for any projective $k$-scheme. Let us then prove that Theorem \ref{teo0} holds for any proper $k$-scheme $X$. We argue by induction on the dimension of $X$, the case where $X$ is $0$-dimensional being immediate.

By Chow's lemma, there exists a proper morphism $f : X' \rightarrow X$, which is an isomorphism above a dense open subset $U$ of $X$, such that $X'$ is projective. Let $Z' = X' \setminus f^{-1}(U)$ and $Z = X \setminus U$. Since $X'$ and $Z'$ are projective, and since $Z$ has dimension strictly smaller than $X$, the preliminary assumption and the induction hypothesis yield that Theorem \ref{teo0} holds for $X',Z'$ and $Z$, hence the existence of collections $(E(X'_{(x')}, -))_{x' \in |X'|}$, $(E(Z'_{(z')},-))_{z' \in |Z'|}$ and $(E(Z_{(z)},-))_{z \in |Z|}$ of triangulated functors satisfying the conclusion of Theorem \ref{teo0} for $X',Z'$ and $Z$ respectively.

 Let $\F$ be an object of $D_c^b(X,\Lambda)$. Then $\varepsilon_{\bk}(X, \F)$ coincides with $\varepsilon_{\bk}(Z, \F) \varepsilon_{\bk}(X', \F) \varepsilon_{\bk}(Z', \F)^{-1}$, and is therefore equal to the determinant of
 \begin{align*}
  \bigoplus_{z \in |Z|} \Ind_{G_z}^{G_k} E(Z_{(z)},\F_{|Z_{(z)}}) \oplus \bigoplus_{x' \in |X'|} \Ind_{G_{x'}}^{G_k} E(X'_{(x')},f^{-1}\F_{|X'_{(x')}}) \oplus \bigoplus_{z' \in |Z'|} \Ind_{G_{z'}}^{G_k} E(Z'_{(z')},f^{-1}\F_{|Z'_{(z')}})[1].
  \end{align*}
This yields the conclusion by setting for any closed point $x$ of $X$,
$$
E(X_{(x)},\G) = E(Z_{(x)},\G_{|Z_{(x)}}) \oplus \bigoplus_{x' \in |f^{-1}(x)|} \Ind_{G_{x'}}^{G_x} \left( E(X'_{(x')},f^{-1}\G_{|X'_{(x')}}) \oplus  E(Z'_{(x')},f^{-1}\G_{|Z'_{(x')}})[1] \right)
$$
if $x$ belongs to $Z$, and
$$
E(X_{(x)},\G) = \bigoplus_{x' \in |f^{-1}(x)|} \Ind_{G_{x'}}^{G_x}  E(X'_{(x')},f^{-1}\G_{|X'_{(x')}}).
$$
otherwise.
\subsection{\label{4.2}} Let $\iota : X \rightarrow \mathbb{P}_k^d$ be a closed immersion into a projective space $\mathbb{P}_k^d$, for some integer $d$. If Theorem \ref{teo0} holds for $ \mathbb{P}_k^d$, so that we have a collection $(E(\mathbb{P}^d_{k,(x)},-))_{x \in |\mathbb{P}_k^d|}$ of triangulated functors, then Theorem \ref{teo0} also holds for $X$ by setting
$$
E(X_{(x)},\G) = E(\mathbb{P}^d_{k,(\iota(x))},\iota_* \G). 
$$

\subsection{\label{4.3}} Let us prove Theorem \ref{teo0}. By \ref{4.1} and \ref{4.2}, it is sufficient to consider the case where $X = \mathbb{P}(V)$ is the projective space parametrizing hyperplanes in a $k$-vector space $V$. We denote by $d+1$ the dimension of $V$, and we argue by induction on $d$, the case $d=0$ being immediate. We henceforth assume that $d$ is greater than $0$.

Let $e_1,e_2$ be a basis of a $2$-dimensional sub-$k$-vector space $W$ of $V$ and let $X'$ be the flag variety parametrizing pairs $([t_1:t_2],H)$, where $[t_1:t_2]$ is a point of $\mathbb{P}^1_k$ and $H$ is an hyperplane in $V$ containing $t_1 e_1 + t_2 e_2$. Let $f : X' \rightarrow \mathbb{P}^1_k$ and $r : X' \rightarrow \mathbb{P}(V)$ be the first and second projections respectively.

Let $\F$ be an object of $D_c^b(\mathbb{P}(V),\Lambda)$. The morphism $r$ is an isomorphism outside the closed subscheme $\iota : \mathbb{P}(V/W) \rightarrow \mathbb{P}(V)$, and is a $\mathbb{P}^1$-bundle above $\mathbb{P}(V/W)$, hence
\begin{align}
\varepsilon_{\bk}( \mathbb{P}(V), \F)  &= \varepsilon_{\bk}(X', r^{-1}\F)  \varepsilon_{\bk}(\mathbb{P}(V/W), \iota^{-1}\F(-1))^{-1} \\
&=\varepsilon_{\bk}( \mathbb{P}^1_k, Rf_* r^{-1} \F)   \varepsilon_{\bk}(\mathbb{P}(V/W), \iota^{-1}\F(-1))^{-1}. \label{eq67}
\end{align}
 By \ref{4.0}, the homomorphism $\varepsilon_{\bk}( \mathbb{P}^1_k, Rf_* r^{-1} \F) $ is the determinant of
$$
\mathrm{CFT}_{t^{-1}} R \Psi_{\id}^{\infty} \G[1](-1) \oplus \mathrm{CFT}_{-t^{-1}}  R \Psi_{\id}^{\infty} \G[1](-1) \oplus \bigoplus_{s \in |\mathbb{P}^1_k|} \Ind_{G_s}^{G_k} \Art_{\pi_s} R\Phi_{\id}^s(\G),
$$
where $\G = Rf_* r^{-1} \F$. Moreover, by \ref{isom6} we have a functorial isomorphism
$$
 \Art_{\pi_s} R\Phi_{\id}^s(\G) \simeq R\Gamma(X'_{\overline{s}},\Art_{\pi_s} R\Phi_{f}^s( r^{-1} \F)).
$$
We similarly have
\begin{align*}
\mathrm{CFT}_{t^{-1}} R \Psi_{\id}^{\infty} \G(-1) &\cong  R\Gamma(X'_{\overline{s}},\mathrm{CFT}_{t^{-1}} R\Psi_{f}^{\infty}( r^{-1} \F(-1))) \\
\mathrm{CFT}_{-t^{-1}} R \Psi_{\id}^{\infty} \G &\cong  R\Gamma(X'_{\overline{s}},\mathrm{CFT}_{-t^{-1}} R\Psi_{f}^{\infty}( r^{-1} \F))
\end{align*}

For each closed point $s = [t_1,t_2]$ of $\mathbb{P}^1_k$, the fiber $X'_s$ is the projective space associated to the quotient of $V \otimes_k k(s)$ by the $k(s)$-line $k(s)(t_1 e_1 + t_2 e_2)$. In particular, the induction hypothesis applies to $X'_s$, and we therefore have a collection $(E(X'_{s,(x')},-))_{x' \in |X'_s|}$ of triangulated functors satisfying the conclusion of Theorem \ref{teo0} for the $s$-scheme $X'_s$. Thus $\varepsilon_{\bk}( \mathbb{P}^1_k, Rf_* r^{-1} \F) $ is the determinant of
\begin{align*}
&\bigoplus_{x' \in |X'_{\infty}|}   \Ind_{G_{x'}}^{G_k}  E(X'_{\infty,(x')}, \mathrm{CFT}_{t^{-1}} R\Psi_{f}^{\infty}( r^{-1} \F(-1))) \\
 \oplus &\bigoplus_{x' \in |X'_{\infty}|}   \Ind_{G_{x'}}^{G_k}  E(X'_{\infty,(x')}, \mathrm{CFT}_{-t^{-1}} R\Psi_{f}^{\infty}( r^{-1} \F)) \\
  \oplus  &\bigoplus_{x' \in |X'|} \Ind_{G_{x'}}^{G_k} E(X'_{f(x'),(x')},\Art_{\pi_{f(x')}} R\Phi_{f}^{f(x')}( r^{-1} \F)).
\end{align*}

One should note that $r^{-1} \F$ is universally locally acyclic with respect to $f$ above a dense open subscheme of $\mathbb{P}^1_k$ by (\cite{SGA412}, [Th. finitude] Th. 2.13), hence $ R\Phi_{f}^s( r^{-1} \F)$ vanishes for all but finitely many closed points $s$ of $\mathbb{P}^1_k$. This and the induction hypothesis ensure that only finitely many terms contribute in the above sum. The conclusion then follows from \ref{eq67} and from the induction hypothesis applied to $\mathbb{P}(V/W)$; in particular, the twist formula \ref{teo0}(2) follows from \ref{tensor} and the induction hypothesis.

\section{Proof of Theorem \ref{teo0.2} \label{prf2} }

In this section, we prove Theorem \ref{teo0.2prime} below, of which Theorem \ref{teo0.2} is an immediate consequence. 

\begin{teo}\label{teo0.2prime} Let $\Lambda$ be either $\overline{\mathbb{F}}_{\ell}$ or $\overline{\mathbb{Q}}_{\ell}$. Let $S$ be a henselian trait of equicharacteristic $p$, with closed point $s$ such that $k(s)$ is perfect, and let $\omega$ be a nonzero meromorphic differential $1$-form on $S$. Let $f : X \rightarrow S$ be a proper morphism. Then there exists a collection
$$
(E_{f,\omega,x} )_{x \in |X_s|},
$$ 
of triangulated functors $E_{f,\omega,x} : D_c^b(X_{(x)}, \Lambda) \rightarrow D_c^b(x,\Lambda)$, indexed by the set of closed points of the special fiber $X_s$, such that:
\begin{enumerate}
\item for any object $\F$ of $D_c^b(X,\Lambda)$, we have $E_{f,\omega,x}( \F_{|X_{(x)}} ) \simeq 0$ for all but finitely many closed points $x$ of $X$;
\item for any closed point $x$ of $X$, any object $\F$ of $D_c^b(X_{(x)},\Lambda)$, and any object $\G$ of $D_c^b(x,\Lambda)$, we have an isomorphism
$$
E_{f,\omega,x}( \F \otimes \mathrm{sp}^{-1} \G) \simeq E_{f,\omega,x}( \F) \otimes \G,
$$
where $\mathrm{sp} : X_{(x)} \rightarrow x$ is the canonical specialization morphism, and this isomorphism is functorial in $\F$ and $\G$;
\item for any object $\F$ of $D_c^b(X,\Lambda)$, we have
\begin{align*}
\varepsilon_{\bs}(S, Rf_* \F, \omega) &= \det \left( \bigoplus_{x \in |X_s|} \Ind_{G_x}^{G_s}E_{f,\omega,x}(\F_{|X_{(x)}} ) \right), \\
a(S, Rf_* \F, \omega) &= \rk \left( \bigoplus_{x \in |X_s|} \Ind_{G_x}^{G_s} E_{f,\omega,x}(\F_{|X_{(x)}} ) \right),
\end{align*}
where $\bs$ is a separable closure of $s$ with Galois group $G_s$, and where $a(S,-, \omega) $ is the Artin conductor, cf. $($ \cite{G19}, $7.2)$.
\end{enumerate}
\end{teo}

Let $\pi$ be a uniformizer of $S$, and let us write $\omega = \alpha d \pi$ for some nonzero element $\alpha$ of the function fields $k(\eta)$ of $S$. We denote by $v(\alpha)$ the valuation of $\alpha$ in the discretely valued field $k(\eta)$. For any object $\G$ of $D_c^b(S,\Lambda)$, we have
\begin{align*}
\varepsilon_{\bs}(S, \G, \omega) &= \varepsilon_{\bs}(S, \G, d \pi) \langle \chi_{\det(\G_{\eta})} \rangle(\alpha) \chi_{\cyc}^{-v(\alpha) \rk(\F_{\eta})} , \\
a(S, \G, \omega) &=  a(S, \G)+ v(\alpha) \rk(\G_{\eta}) ,
\end{align*}
with notation as in $($ \cite{G19}, $9.2)$ or Proposition \ref{CFT}. It is therefore enough to produce families $(E^1_{f,x} )_{x \in |X_s|}$ and $(E^2_{f,x} )_{x \in |X_s|}$ of triangulated functors as in Theorem \ref{teo0.2prime} (1),(2), such that
\begin{align*}
\varepsilon_{\bs}(S, Rf_* \F, d \pi) &= \det \left( \bigoplus_{x \in |X_s|} \Ind_{G_x}^{G_s}E^1_{f,x}(\F_{|X_{(x)}} ) \right), \\
a(S, Rf_* \F) &= \rk \left( \bigoplus_{x \in |X_s|} \Ind_{G_x}^{G_s} E^1_{f,x}(\F_{|X_{(x)}} ) \right),
\end{align*}
and
\begin{align*}
\langle \chi_{\det(\G_{\eta})} \rangle(\alpha) \chi_{\cyc}^{-v(\alpha) \rk(\F_{\eta})}  &= \det \left( \bigoplus_{x \in |X_s|} \Ind_{G_x}^{G_s}E^2_{f,x}(\F_{|X_{(x)}} ) \right), \\
v(\alpha) \rk(\G_{\eta})  &= \rk \left( \bigoplus_{x \in |X_s|} \Ind_{G_x}^{G_s} E^2_{f,x}(\F_{|X_{(x)}} ) \right).
\end{align*}
In order to construct $(E^1_{f,x} )_{x \in |X_s|}$, we first notice that Corollary \ref{art10} and Proposition \ref{isom6} yields that $\varepsilon_{\bs}(S, Rf_* \F, d \pi)$ and $a(S, Rf_* \F)$ are respectively the determinant and rank of the complex
$$
\Art_{\pi} R\Phi_{\id}^s(Rf_* \F) \cong R \Gamma(X_{\bs}, \Art_{\pi} R\Phi_{f}^s(\F)). 
$$
If $(E(X_{s,(x)}, - ) )_{x \in |X_s|}$ is as in Theorem \ref{teo0} applied to the proper $s$-scheme $X_s$, then we can set
$$
E^1_{f,x}(\F) = E(X_{s,(x)},  \Art_{\pi} R\Phi_{f}^s(\F) ).
$$ 
For $(E^2_{f,x} )_{x \in |X_s|}$, we can assume by additivity that $\alpha$ is a uniformizer on $S$, in which case Propositions \ref{CFT} ensures that it is enough to set
$$
E^2_{f,x}(\F) = E(X_{s,(x)},  \mathrm{CFT}_{\alpha} R\Psi_{f}^s(\F)(-1) ).
$$ 

\begin{rema} Let $X$ be a smooth projective curve of a field $k$ of positive characteristic $p$ and let $\F$ be an object of $D_c^b(X,\Lambda)$. The product formula for curves expresses the determinant of $R \Gamma(X_{\bk},\F)$ as a product of finitely many local $\varepsilon$-factors, which are invariants attached to Galois representations of $1$-dimensional local fields. If $k$ is itself a $1$-dimensional local field, then the method described in this section decomposes the local $\varepsilon$-factor of $R \Gamma(X_{\bk},\F)$ as a product of finitely many contributions, namely $2$-dimensional local $\varepsilon$-factors, which are invariants attached to Galois representations of $2$-dimensional local fields. By iterating this procedure, one can attach an $n$-dimensional local $\varepsilon$-factor to Galois representations of $n$-dimensional local fields, such that if $k$ is an $n$-dimensional local field then the $n$-dimensional local $\varepsilon$-factor of $R \Gamma(X_{\bk},\F)$ factors as the product of finitely many $(n+1)$-dimensional local $\varepsilon$-factors.
\end{rema}

\bibliographystyle{amsalpha}

\begin{thebibliography}{2}
%
%
%
%
%
%
%
%
%
%
\bibitem[De73]{De73} P. Deligne, ``Les constantes des \'equations fonctionnelles des fonctions $L$'', in ``Modular Functions of One Variable II'', Lecture Notes in Mathematics 349, Springer-Verlag, 1973.
%
%
%
\bibitem[Dw56]{Dw} B. Dwork, ``On the Artin root number'', Amer. J. Math. 78, 1956, pp.444-472.
%
%
%
%
%
\bibitem[SGA4]{SGA4}
A. Grothendieck, \emph{S\'eminaire de G\'eom\'etrie Alg\'ebrique du Bois Marie - 1963-64 - Th\'eorie des topos et cohomologie \'etale des sch\'emas - (SGA 4)}, Springer-Verlag, LNM 269/270/305, 1972/3. 

\bibitem[SGA $4 \frac{1}{2}$]{SGA412} P. Deligne, \emph{S\'eminaire de G\'eom\'etrie Alg\'ebrique du Bois Marie - Cohomologie \'etale - (SGA 4 1/2)}, Springer-Verlag, LNM 569, 1977.


%

\bibitem[SGA7]{SGA7} P. Deligne and N. Katz, \emph{S\'eminaire de G\'eom\'etrie Alg\'ebrique du Bois Marie - 1967-69 - Groupes de monodromie en g\'eom\'etrie alg\'ebrique - (SGA 7)}, vol. 2, Springer-Verlag, LNM 340, 1973.

%

\bibitem[Gu19]{G19} Q. Guignard, ``Geometric local epsilon factors'', arxiv:1902.06523 (v3).

\bibitem[Ka86]{katz} N. M. Katz, ``Local-to-global extensions of representations of fundamental groups'', Annales de l'institut Fourier, tome 36, No. 4, 69-106, 1986.
%
%
%
%

\bibitem[Ill81]{Ill81} Luc Illusie, ``Th\'eorie de Brauer et caract\'eristique d'Euler-Poincar\'e d'apr\`es P. Deligne'', Ast\'erisque 82-83 (1981), pp. 161-172.

\bibitem[KS]{KS} K. Kato and S. Saito, ``Unramified class field theory of arithmetical surfaces'', Annals of Mathematics 118 (1983), pp. 241-275.

%
%
\bibitem[La83]{La83} G. Laumon, ``Vanishing cycles over a base of dimension $\geq 1$'', Algebraic geometry (Tokyo/Kyoto, 1982), Lecture Notes in Mathematics, vol. 1016, Springer, Berlin, 1983, pp. 143-150. 

\bibitem[La87]{La87} G. Laumon, ``Transformation de Fourier, constantes d'\'equations fonctionnelles et conjecture de Weil'', Publications math\'ematiques de l'I.H.E.S., tome 65 (1987), pp. 131-210.
%
%
%
%
%
\bibitem[Sa17]{Sai} T. Saito, ``The characteristic cycle and the singular support of a constructible sheaf'', Inventiones mathematicae 207(2) (2017), pp. 597-695.

%
%
%
%
\bibitem[Se98]{S} J.-P. Serre, ``Repr\'esentations lin\'eaires des groupes finis'', Hermann, M\'ethodes, 1998
%
%
%
%
\bibitem[ST68]{ST} J.-P. Serre, J. Tate, ``Good Reduction of Abelian Varieties'', Annals of Mathematics Second Series, Vol. 88, No. 3, 492-517, 1968.

\bibitem[Stacks]{stacks-project}
The {Stacks Project Authors}, \emph{Stacks {P}roject},
 \url{http://stacks.math.columbia.edu}, 2019.
%
%
\bibitem[UYZ]{UYZ} N. Umezaki, E. Yang, Y. Zhao, ``Characteristic class and the epsilon factor of an \'etale sheaf'', arxiv:1701.02841.


\bibitem[Ya1]{Ya1} S. Yasuda, ``Local Constants in Torsion Rings'', Journal of Mathematical Sciences (University of Tokyo), Vol. 16, No. 2, 125-197, 2009.

\bibitem[Ya2]{Ya2} S. Yasuda, ``The Product Formula for Local Constants in Torsion Rings'', Vol. 16, No. 2, 199-230, 2009.

\bibitem[Ya3]{Ya3} S. Yasuda, ``Local $\varepsilon_0$-characters in torsion rings'', Journal de Th\'eorie des Nombres
de Bordeaux, Vol. 19, 763-797, 2007.

\end{thebibliography}

\end{document}